\documentclass{amsart}
\usepackage{pb-diagram,lamsarrow,pb-lams,latexsym,graphics}

\makeatletter
\def\l@subsection{\@tocline{2}{0pt}{3.3pc}{5pc}{}}
\makeatother

\makeatletter
\def\subsection{\@startsection{subsection}{2}%
  \z@{.7\linespacing\@plus\linespacing}{-.5em}%
  {\normalfont\bfseries}}
\makeatother

\newcommand{\vspacea}{\vspace{3pt}}
\newcommand{\vspaceb}{\vspace{5pt}}
\newcommand{\vspacec}{\vspace{7pt}}

\theoremstyle{plain}
\newtheorem{Theorem}{Theorem}
\newtheorem{Corollary}[Theorem]{Corollary}
\newtheorem{Lemma}[Theorem]{Lemma}
\newtheorem*{DLemma}{Depth Lemma}
\newtheorem{Proposition}[Theorem]{Proposition}
\newtheorem{Definition}[Theorem]{Definition}

\theoremstyle{remark}
\newtheorem*{Remark}{Remark}
\newtheorem*{Remarks}{Remarks}
\newtheorem*{Notation}{Notation}

\DeclareMathOperator{\Spec}{Spec}
\DeclareMathOperator{\Supp}{Supp}

\DeclareMathOperator{\Hom}{Hom}
\DeclareMathOperator{\Ext}{Ext}
\DeclareMathOperator{\Ker}{Ker}
\DeclareMathOperator{\Coker}{Coker}
\DeclareMathOperator{\Img}{Im}
\DeclareMathOperator{\idmap}{id}
\DeclareMathOperator{\Tor}{Tor}

\DeclareMathOperator{\Ann}{Ann}
\DeclareMathOperator{\depth}{depth}
\DeclareMathOperator{\grade}{grade}
\DeclareMathOperator{\gdim}{G-dim}
\DeclareMathOperator{\pd}{pd}
\DeclareMathOperator{\injdim}{inj.dim}

\newcommand{\F}{\mathbf{f}}

\newcommand{\kk}{\mathbf{k}}
\newcommand{\mm}{\mathfrak{m}}

\newcommand{\union}{\cup}

\newcommand{\tens}{\otimes}

\newcommand{\isom}{\cong}

\newcommand{\eqdef}{\overset{\textup{def}}{=}}
\newcommand{\proje}{\approx}
\newcommand{\lloc}{<_{\textup{loc}}}
\newcommand{\vertisom}{\rotatebox{270}{$\isom$}}

\begin{document}



\title[Gorenstein dimension]{Gorenstein dimension of modules}

\author{Vladimir Ma\c{s}ek}
\address{Department of Mathematics, Box 1146, Washington University,
        St. Louis, MO 63130}
\email{vmasek@math.wustl.edu}


\maketitle




\tableofcontents

\noindent {\bf Notation.}

\vspacea

\begin{tabbing}
99\=9999999999999999999999\=9999999999999999999999999999\kill
  \>$R$ \> ring (always commutative and Noetherian)                 \\
  \>$(R,\mm,\kk)$ \> local ring with maximal ideal $\mm$ and $\kk=R/{\mm}$ \\
  \>$L,M,N,\dots$     \> $R$-modules (always finitely generated)    \\
  \>$M^*$             \> $\Hom_R(M,R)$, the dual of $M$   \\
  \>$D(M)$   \> the Auslander dual of $M$ (Definition \ref{def2}) \\
  \>$\sigma_M:M\arrow{e} M^{**}$ \> the natural evaluation map; \\
  \>  \hspace{28pt}  $K_M = \Ker(\sigma_M), C_M = \Coker(\sigma_M)$  \\
  \>$\gdim_R(M),\gdim(M)$ \> Gorenstein dimension of $M$
                (Definition \ref{def16})  \\
  \>$\gdim(M) \lloc \infty$ \> $M$ has locally finite Gorenstein dimension \\
  \>                        \> (Remark \#2 after Theorem \ref{ABF})   \\
  \>$\pd_R(M),\pd(M)$     \>  projective dimension of $M$ \\
  \>$\proje$     \>  projective equivalence (Definition \ref{def3}) \\
\end{tabbing}

\vspaceb


\section{Introduction}

In these expository notes I discuss several concepts and results in the theory
of modules over commutative rings, revolving around the Gorenstein dimension
of modules. Some of the related notions are the \emph{Auslander dual},
\emph{$k$-torsionless modules}, and \emph{$k^{\text{th}}$ syzygies}.

Essentially everything in these notes can be found, in one form or another, in
the memoir ``Stable module theory'' by M.~Auslander and M.~Bridger (Mem. A.M.S.,
no. 94, 1969). The only difference is in presentation. In the Auslander--Bridger
memoir many of the results are proved in the most general setting, e.g. over
possibly non-commutative, non-Noetherian rings. The techniques used are quite
abstract and unfamiliar to many commutative algebraists. Much space is devoted
to the theory of satellites of functors which are exact only in the middle, etc.
While such a degree of generality has many advantages, it does make the memoir
difficult to read for the non-expert.
My goal in writing these notes was to develop the theory in the context of
commutative Noetherian rings, and to show that, in this important special case,
the theory is fairly elementary and easy to build. As a practical matter, then,
I wrote the notes using Matsumura's ``Commutative ring theory'' as the only
prerequisite; and indeed, my hope is that these notes can be read just like
an extra chapter in Matsumura's book. Still, some of the proofs given here are
mere adaptations and simplification of those in \cite{ab}. Incidentally, in
\S2 I fix a mistake in the proof of the analogue of the Auslander--Buchsbaum
formula for G-dimension given in \cite{ab}.

\vspacec

{\bf Throughout these notes all rings are commutative and Noetherian and
all modules are finitely generated.}

\vspacec

\noindent {\bf Motivation.}
Here are some of the reasons why anyone would want to study the theory of
Gorenstein dimension of modules:
  \begin{itemize}
    \addtolength{\itemsep}{3pt}
    \item Modules of finite projective dimension also have finite
      Gorenstein dimension. On the other hand, \emph{all} modules over a
      Gorenstein ring have locally finite Gorenstein dimension. There are
      many results in commutative algebra which begin with a hypothesis
      like, ``assume that either $\pd(M)<\infty$ or that the ring $R$ is
      Gorenstein.'' And indeed, upon inspection, it turns out that the
      results hold for modules of (locally) finite Gorenstein dimension
      over any ring. The Gorenstein dimension provides a natural
      unifying language for such results.
    \item Gorenstein dimension has many of the good properties of
      projective dimension. (In particular, the Auslander--Buchsbaum
      formula holds with Gorenstein dimension instead of projective
      dimension, and a local ring $(R,\mm,\kk)$ is Gorenstein if and only
      if $\gdim(\kk)<\infty$. Compare with: ``$R$ is regular if and only if
      $\pd(\kk)<\infty$.'') For this reason, many results proved initially
      over regular rings hold, in fact, over Gorenstein rings. One such
      example is the existence of Evans--Griffith presentations \cite{mas}.
    \item Some of the most useful characterizations of syzygies hold for
      modules of finite Gorenstein dimension (Theorem~\ref{thm40} in \S 3).
    \item In \cite{af}, Avramov and Foxby defined local ring homomorphisms
      of finite Gorenstein dimension and studied their properties. In
      particular, they defined a dualizing complex for such homomorphisms,
      very similar to the dualizing complex of a local ring.
  \end{itemize}

\noindent {\bf Acknowledgement.}
In preparing and writing up these notes, I benefited from advice, encouragement
and many useful discussions with L.~Avramov, E.~G.~Evans, P.~Griffith, N.~Mohan
Kumar, and D.~Popescu. The treatment of several topics here owes much to their
insights.


\section{The Auslander dual and $k$-torsionless modules}

\begin{Definition} \label{def1}
    A module $M$ is \emph{torsionless}, resp. \emph{reflexive}, if the
    natural evaluation map $\sigma_M:M\arrow{e} M^{**}$ is injective, resp.
    an isomorphism. That is, $M$ is torsionless if $K_M=0$, and reflexive
    if $K_M=C_M=0$, where $K_M=\Ker(\sigma_M)$ and $C_M=\Coker(\sigma_M)$.
    \emph{(See Remark after Definition~\ref{def35} in \S 3 for a short
    comment on terminology.)}
\end{Definition}

Auslander interpreted $K_M$ and $C_M$ cohomologically, in terms of the
``Auslander dual'', $D(M)$, defined below. This interpretation explains the
functorial properties of $K_M$ and $C_M$ and the good behavior of torsionless
and reflexive modules, and can be used to define ``$k$-torsionless'' for all
$k \geq 0$; see Definition \ref{def7} later in this section. ($1$-torsionless
is the same as torsionless, and $2$-torsionless is the same as reflexive.
\emph{All} modules are $0$-torsionless.)

\begin{Definition} \label{def2}
    Let $M$ be any module (finitely generated as always), and let
      \begin{equation}
        P_1 \arrow{e,t}{u} P_0 \arrow{e,t}{f} M \arrow{e} 0 \tag{$\pi$}
      \end{equation}
    be a (finitely generated) projective presentation of $M$.

    The \emph{Auslander dual}, $D(M)$, of $M$, is defined as
      \[
        D(M) = \Coker(u^*:P_0^* \arrow{e} P_1^*);
      \]
    in other words, dualizing $(\pi)$ we get an exact sequence
      \begin{equation}
        0 \arrow{e} M^* \arrow{e,t}{f^*} P_0^* \arrow{e,t}{u^*}
          P_1^* \arrow{e} D(M) \arrow{e} 0.  \tag{$\pi^*$}
      \end{equation}
\end{Definition}

$D(M)$ is sometimes called the \emph{Auslander transpose} of $M$; see, for
example, \cite[p. 648]{eis}.

\vspaceb

Clearly, $D(M)$ \emph{depends} on which projective presentation ($\pi$) is
used in the definition. Until the end of the proof of uniqueness of $D(M)$
(up to projective equivalence) in Proposition \ref{prop4} below,
we will denote $\Coker(u^*)$ by $D_{\pi}(M)$.
    
\begin{Definition} \label{def3}
    Two modules $M$ and $N$ are \emph{projectively equivalent} if $\exists
    \, P,Q$ projective with $M \oplus P \cong N \oplus Q$. \emph{Notation:}
    \mbox{$M \proje N$.}
\end{Definition}

$\proje$ is an equivalence relation on the class of (finitely generated)
$R$-modules.

\begin{Proposition} \label{prop4}
    If $(\pi)$ as above and
      \begin{equation}
        Q_1 \arrow{e,t}{v} Q_0 \arrow{e,t}{g} M \arrow{e} 0 \tag{$\rho$}
      \end{equation}
    are two projective presentations of $M$, then
        \[
            D_{\pi}(M) \proje D_{\rho}(M).
        \]
\end{Proposition}
 
\begin{proof}
  We say that ($\pi$) \emph{strictly dominates} ($\rho$) if there are linear
  maps $\phi_i : P_i \arrow{e} Q_i \; (i=0,1)$ such that
      \begin{enumerate}
        \item[(a)] $\phi_i$ is surjective, $i=0,1$;
        \item[(b)] $\phi_0$ is a lifting of $\idmap_M$ and $\phi_1$
               is a lifting of $\phi_0$; i.e., the diagram
  \begin{equation*}
    \begin{diagram}          
        \node{P_1} \arrow{e,t}{u} \arrow{s,r}{\phi_1}
        \node{P_0} \arrow{e,t}{f} \arrow{s,r}{\phi_0}
        \node{M} \arrow{e} \arrow{s,=} \node{0}     \\
        \node{Q_1} \arrow{e,b}{v} \node{Q_0} \arrow{e,b}{g}
        \node{M} \arrow{e} \node{0}
    \end{diagram}            
  \end{equation*}
             commutes;
        \item[(c)] The map $\bar{u} : K_1 \arrow{e} K_0$ induced by $u$
        is surjective, where $K_i = \Ker(\phi_i) \; (i=0,1)$.
      \end{enumerate}
  In other words, we should have a commutative diagram with exact rows and
  columns:
  \begin{equation*} \label{diag:**}
    \begin{diagram}          
        \node{0} \arrow{s} \node{0} \arrow{s}   \\
        \node{K_1} \arrow{e,t}{\bar{u}} \arrow{s}
        \node{K_0} \arrow{e} \arrow{s} \node{0}\\
        \node{P_1} \arrow{e,t}{u} \arrow{s,r}{\phi_1}
        \node{P_0} \arrow{e,t}{f} \arrow{s,r}{\phi_0}
        \node{M} \arrow{e} \arrow{s,=} \node{0}     \\
        \node{Q_1} \arrow{e,b}{v} \arrow{s}
        \node{Q_0} \arrow{e,b}{g} \arrow{s}
        \node{M} \arrow{e} \node{0}     \\
        \node{0}  \node{0} 
    \end{diagram}            \tag*{$(\ast\ast)$}
  \end{equation*}
  We prove the Proposition in two steps. First we show that $D_{\pi}(M) \proje
  D_{\rho}(M)$ under the extra hypothesis that ($\pi$) strictly dominates
  ($\rho$). Then we show that, given any two projective presentations ($\pi$)
  and ($\rho$), there is a third one, ($\sigma$), which strictly dominates
  both ($\pi$) and ($\rho$).

  \vspaceb

  For the time being, assume that ($\pi$) strictly dominates ($\rho$). Then we
  have the diagram \ref{diag:**}. Since the $Q_i$ are projective, the columns
  are split-exact and the $K_i$ are projective. Therefore $\bar{u}$ splits.

  Dualizing \ref{diag:**} we get a commutative diagram with exact rows and
  columns:
  \begin{equation*}
    \begin{diagram}          
        \node[3]{0} \arrow{s} \node{0} \arrow{s} \node{0}\arrow{s}  \\
        \node{0} \arrow{e} \node{M^*} \arrow{e} \arrow{s,=}
        \node{Q_0^*} \arrow{e} \arrow{s} \node{Q_1^*} \arrow{e}
        \arrow{s} \node{D_{\rho}(M)} \arrow{e} \arrow{s} \node{0} \\
        \node{0} \arrow{e} \node{M^*} \arrow{e} \node{P_0^*} \arrow{e}
        \arrow{s} \node{P_1^*} \arrow{e} \arrow{s}  \node{D_{\pi}(M)}
        \arrow{e} \arrow{s} \node{0}            \\
        \node[2]{0} \arrow{e} \node{K_0^*} \arrow{e,t}{\bar{u}^*}
        \arrow{s} \node{K_1^*} \arrow{e} \arrow{s} \node{K} \arrow{e}
        \arrow{s} \node{0}              \\
        \node[3]{0}  \node{0}  \node{0}
    \end{diagram}            
  \end{equation*}
  where $K \eqdef \Coker(\bar{u}^*)$. \ \ $0 \to D_{\rho}(M) \to
  D_{\pi}(M) \to K \to 0$ is exact by the Snake Lemma, and $K$ is projective,
  because $\bar{u}^*$ is a split injective map of projective modules. Thus
  $D_{\pi}(M) \cong D_{\rho}(M) \oplus K$, and $D_{\pi}(M) \proje D_{\rho}(M)$.

  \vspaceb

  Now let ($\pi$), ($\rho$) be any two projective presentations of $M$. We
  construct a new presentation, ($\sigma$), which strictly dominates both
  ($\pi$) and ($\rho$).

  Let $h:P_0 \oplus Q_0 \arrow{e} M,\,h(p_0,q_0)=f(p_0)+g(q_0)$.
  Clearly, $h$ is surjective. Let $\alpha : E \arrow{e} P_0 \oplus Q_0$ be a
  linear map from a (finitely generated) projective module $E$ onto
  $\Ker(h) \subseteq P_0 \oplus Q_0$; that is, consider an exact sequence
  \[
      E \arrow{e,t}{\alpha} P_0 \oplus Q_0 \arrow{e,t}{h} M \arrow{e} 0.
  \]
  Extend the projective presentations ($\pi$) and ($\rho$) of $M$ to the left:
  \[
      P_2 \arrow{e,t}{u'} P_1 \arrow{e,t}{u} P_0 \arrow{e,t}{f} M \arrow{e} 0
  \]
  and
  \[
      Q_2 \arrow{e,t}{v'} Q_1 \arrow{e,t}{v} Q_0 \arrow{e,t}{g} M \arrow{e} 0.
  \]
  Define the projective presentation ($\sigma$) of $M$ as follows:
    \begin{equation*}
       E \oplus P_2 \oplus Q_2 \arrow{e,t}{w} P_0 \oplus Q_0 \arrow{e,t}{h}
       M \arrow{e} 0,       \tag{$\sigma$}
    \end{equation*}
  where $w(e,p_2,q_2) \eqdef \alpha(e)$. \ \ \ ($\sigma$) is clearly exact.

  I claim that ($\sigma$) strictly dominates both ($\pi$) and ($\rho$). Since
  the construction of ($\sigma$) is symmetric with respect to ($\pi$) and
  ($\rho$), I will show only that ($\sigma$) strictly dominates ($\rho$).

  Lift $\idmap_M$ to $\phi_0:P_0\to Q_0$; that is, $g\,\phi_0=f$. Define
  $\chi_0 : P_0 \oplus Q_0 \to Q_0, \, \chi_0(p_0,q_0)=\phi_0(p_0)+q_0$.
  $\chi_0$ is clearly surjective, and $g \chi_0 = h$.

  Lift $\chi_0$ to $\delta : E \to Q_1$; i.e., construct a commutative
  diagram
  \begin{equation*}
    \begin{diagram}          
        \node{E} \arrow{e,t}{\alpha} \arrow{s,r}{\delta}
        \node{P_0 \oplus Q_0} \arrow{e,t}{h} \arrow{s,r}{\chi_0}
        \node{M} \arrow{e} \arrow{s,=} \node{0\,}       \\
        \node{Q_1} \arrow{e,b}{v} \node{Q_0} \arrow{e,b}{g}
        \node{M} \arrow{e} \node{0.}
    \end{diagram}            
  \end{equation*}

  Finally, define $\chi_1 : E \oplus P_2 \oplus Q_2 \to Q_1, \,
  \chi_1(e,p_2,q_2) = \delta(e)+v'(q_2)$. Then $v \chi_1(e,p_2,q_2) =
  v \delta(e) + v v'(q_2) = \chi_0 \alpha(e) = \chi_0 w(e,p_2,q_2)$; 
  that is, ($\sigma$) and ($\rho$) sit in a commutative diagram:
  \begin{equation*}
    \begin{diagram}          
        \node{E \oplus P_2 \oplus Q_2} \arrow{e,t}{w}
                        \arrow{s,r}{\chi_1}
        \node{P_0 \oplus Q_0} \arrow{e,t}{h} \arrow{s,r}{\chi_0}
        \node{M} \arrow{e} \arrow{s,=} \node{0\,}       \\
        \node{Q_1} \arrow{e,b}{v} \node{Q_0} \arrow{e,b}{g}
        \node{M} \arrow{e} \node{0.}
    \end{diagram}            
  \end{equation*}

  $\chi_1$ is surjective. Indeed, take any $q_1\in Q_1$. Then $v(q_1)\in Q_0$,
  and $h(0,v(q_1))=f(0)+g v(q_1)=0$; so $(0,v(q_1)) \in \Ker(h) = \Img(\alpha)$.
  Take $e \in E$ with $\alpha(e)=(0,v(q_1))$. Then $v\delta(e) =
  \chi_0\alpha(e) = \chi_0(0,v(q_1)) = v(q_1)$; thus $q_1-\delta(e) \in \Ker(v)
  = \Img(v')$, i.e. $q_1 = \delta(e) + v'(q_2)=\chi_1(e,0,q_2)$ for a suitable
  $q_2 \in Q_2$, and we showed that $q_1 \in \Img(\chi_1)$.

  Let $K_i \eqdef \Ker(\chi_i),\,i=0,1$; to finish the proof, we show that the
  map $\bar{w} : K_1 \to K_0$, induced by $w$, is surjective.

  Let $(p_0,q_0) \in K_0 = \Ker(\chi_0) \subseteq \Ker(h) = \Img(\alpha)$.
  Take $e \in E$ such that $\alpha(e) = (p_0,q_0)$.
  $v\delta(e) = \chi_0 \alpha(e) = 0 \implies \delta(e) \in\Ker(v)=\Img(v')$;
  take $q_2 \in Q_2$ such that $\delta(e) = v'(q_2)$. Then $(e,0,-q_2) \in
  \Ker(\chi_1) = K_1$, and $w(e,0,-q_2) = \alpha(e) = (p_0,q_0)$. We showed
  that $w$ takes $K_1$ onto $K_0$, as required.
\end{proof}

\begin{Remarks}
  (1)  When $M$ is given, $D(M)$ is defined only up to projective equivalence.
       However, we will work with $D(M)$ loosely, as though it were an
       $R$-module. We will be careful to specify, when necessary, that a
       particular representative is being used. In many instances, e.g. in
       definitions depending only on vanishing of various $\Ext^i(D(M),R),
       i\geq 1$, the distinction is irrelevant.

  (2)  If $M \proje 0$ (that is, if $M$ is projective), then $D(M) \proje 0$.
       Also, $D(M_1 \oplus M_2) \proje D(M_1) \oplus D(M_2)$. Therefore
       $M \proje N \implies D(M) \proje D(N)$.

  (3)  If we use ($\pi$) to define $D(M)$ (or, rather, $D_{\pi}(M)$), then
       ($\pi^*$) is a projective presentation of $D_{\pi}(M)$, and 
       $D_{\pi^*}(D_{\pi}(M)) = \Coker(u^{**}) \isom M$ (because $P_i \isom
       P_i^{**}$, and $u^{**}$ is identified canonically with $u$). Dropping
       the dependence on ($\pi$) and ($\pi^*$), we can write $D(D(M))\proje M$.

  (4)  If $R$ is local, then we can give a more intrinsic definition of $D(M)$,
       requiring that ($\pi$) be a \emph{minimal} projective presentation of
       $M$. Then $D(M)$ is well-defined up to isomorphism (rather than
       projective equivalence). Note, however, that even then we will not have
       $D(D(M)) \isom M$ in general, since the dual of a minimal projective
       presentation is not necessarily minimal.

  (5)  The Auslander dual commutes with base change. That is, if $R \arrow{e}
       R'$ is any homomorphism of (commutative Noetherian) rings and $M$ is a 
       (finitely generated) $R$-module, then $D_{R'}(M\tens_R R') \proje D_R(M)
       \tens_R R'$ (we use the ring as a subscript to $D$ for emphasis).
       For example, if $S$ is a multiplicative system in $R$, then
       $D_{S^{-1}R}(S^{-1}M) \proje S^{-1}D_R(M)$. In particular, $D_{R_P}(M_P)
       \proje D_R(M)_P,\, \forall P \in \Spec(R)$.
\end{Remarks}

\vspacec

Next we study the relationship between $D(M)$ and $K_M$ and $C_M$:

\begin{Proposition} \label{prop5}
  Let $\sigma_M : M \arrow{e} M^{**}$ be the natural map, with kernel $K_M$
  and cokernel $C_M$. 

  Then we have natural isomorphisms
  \[
      K_M \isom \Ext^1(D(M),R) \quad \text{\emph{and}} \quad C_M
        \isom \Ext^2(D(M),R).
  \]

  Moreover, $\Ext^i(D(M),R) \isom \Ext^{i-2}(M^*,R), \, \forall i \geq 3$.
\end{Proposition}

\noindent
Note that the $\Ext^i(D(M),R)$ do not depend on which particular $D(M)$
is being used.

\begin{proof}
  Consider the projective presentation ($\pi$) of $M$, as before. Dualizing
  ($\pi$) we get ($\pi^*$). Split ($\pi^*$) into short exact sequences:
  \begin{equation*}
    0 \arrow{e} M^* \arrow{e,t}{f^*} P_0^*\arrow{e,t}{\beta_0} N \arrow{e} 0
    \tag{$\pi_0^*$}
  \end{equation*}
  and
  \begin{equation*}
    0 \arrow{e} N \arrow{e,t}{\beta_1} P_1^* \arrow{e} D(M) \arrow{e} 0,
                                \tag{$\pi_1^*$}
  \end{equation*}
  where $N \eqdef \Coker(f^*)$ and $\beta_1 \circ \beta_0 = u^*$.

  Dualizing ($\pi_0^*$) we get an exact sequence
  \begin{equation*}
    0 \arrow{e} N^* \arrow{e,t}{\beta_0^*} P_0^{**} \arrow{e,t}{f^{**}} M^{**}
        \arrow{e} \Ext^1(N,R) \arrow{e} 0,      \tag{$\pi_0^{**}$}
  \end{equation*}
  and dualizing ($\pi_1^*$) we get an exact sequence
  \begin{equation*}
    0 \arrow{e} D(M)^* \arrow{e} P_1^{**} \arrow{e,t}{\beta_1^*} N^*
        \arrow{e} \Ext^1(D(M),R) \arrow{e} 0.       \tag{$\pi_1^{**}$}
  \end{equation*}

  Consider the commutative diagram with exact rows:
  \begin{equation*}
    \begin{diagram}
       \node[2]{P_1} \arrow{e,t}{u} \arrow{s,r}{\beta_1^* \circ \sigma_{P_1}}
       \node{P_0} \arrow{e,t}{f} \arrow{s,r}{\sigma_{P_0}} \node{M} \arrow{e}
       \arrow{s,r}{\sigma_M} \node{0}       \\
       \node{0} \arrow{e} \node{N^*} \arrow{e,b}{\beta_0^*} \node{P_0^{**}}
       \arrow{e,b}{f^{**}} \node{M^{**}}
    \end{diagram}
  \end{equation*}
  Since $\sigma_{P_1}$ is an isomorphism, we have $\Coker(\beta_1^* \circ 
  \sigma_{P_1}) = \Coker(\beta_1^*) \isom \Ext^1(D(M),R)$, the isomorphism
  being given by ($\pi_1^{**}$). As $\sigma_{P_0}$ is an isomorphism, the
  Snake Lemma gives $K_M \eqdef \Ker(\sigma_M) \isom \Coker(\beta_1^* \circ
  \sigma_{P_1}) \isom \Ext^1(D(M),R)$.

  On the other hand, since $f$ is surjective and $\sigma_{P_0}$ is an
  isomorphism, we get $\Img(f^{**}) = \Img(\sigma_M)$, and therefore $C_M
  \eqdef \Coker(\sigma_M) = \Coker(f^{**}) \isom \Ext^1(N,R) \isom
  \Ext^2(D(M),R)$ (the last isomorphism comes from ($\pi_1^*$); the one
  before it from ($\pi_0^{**}$)).

  The last statement of the Proposition follows from ($\pi_0^{*}$) and
  ($\pi_1^{*}$).
\end{proof}

\begin{Remark}
  Essentially the same proof shows the existence of an exact sequence of
  functors:
  \[
    0 \to \Ext^1(D(M),\star) \to M \otimes \star
      \arrow{e,t}{\sigma_M^\star}
      \Hom(M^*,\star) \to \Ext^2(D(M),\star) \to 0;
  \]
  that is, for every $R$-module $N$ there is a natural exact sequence
  \[
    0 \to \Ext^1(D(M),N) \to M \otimes N \arrow{e,t}{\sigma_M^N}
      \Hom(M^*,N) \to \Ext^2(D(M),N) \to 0.
  \]
  (The evaluation map $\sigma_M^N$ is defined by $\sigma_M^N(m\otimes n) =
  \phi_{m,n}, \; \phi_{m,n}(u) = u(m) \cdot n, \; \forall u \in M^*$. The
  naturality is with respect to homomorphisms $N \arrow{e} N'$.) We leave the
  proof as an exercise.
\end{Remark}

\begin{Lemma} \label{lem6}
  Let 
  \[
    0 \arrow{e} M' \arrow{e} M \arrow{e} M'' \arrow{e} 0
  \]
  be an exact sequence. Then, for a suitable choice of Auslander duals,
  we have a long exact sequence:
  \[
    0 \to {M''}^* \to M^* \to {M'}^* \to D(M'') \to D(M) \to D(M') \to 0.
  \]
\end{Lemma}

\begin{proof}
  Let
  \[
    P_1' \arrow{e} P_0' \arrow{e} M' \arrow{e} 0
  \]
  and
  \[
    P_1'' \arrow{e} P_0'' \arrow{e} M'' \arrow{e} 0
  \]
  be any projective presentations of $M'$ and $M''$, respectively. Fit these
  in a commutative diagram with exact rows and columns: 
    \begin{equation*}
      \begin{diagram}       
         \node{0} \arrow{s} \node{0} \arrow{s} \node{0} \arrow{s} \\
         \node{P_1'} \arrow{e} \arrow{s} \node{P_0'} \arrow{e} \arrow{s}
           \node{M'} \arrow{e} \arrow{s} \node{0}         \\
         \node{P_1' \oplus P_1''} \arrow{e} \arrow{s} 
         \node{P_0' \oplus P_0''} \arrow{e} \arrow{s} 
         \node{M} \arrow{e} \arrow{s} \node{0}            \\
         \node{P_1''} \arrow{e} \arrow{s} \node{P_0''} \arrow{e} \arrow{s}
           \node{M''} \arrow{e} \arrow{s} \node{0}        \\
         \node{0}  \node{0}  \node{0}
      \end{diagram}       
    \end{equation*}
  
  The first two columns are split exact, so they remain exact after dualizing.
  Dualize the whole diagram, writing the cokernels of the dual rows as
  $D(M'')$, $D(M)$, and $D(M')$; the conclusion follows from the Snake Lemma.
\end{proof}

\begin{Remark}
  Let $C = \Coker(M^* \to {M'}^*)$. $0 \to M' \to M \to M'' \to 0$ is
  \emph{dual exact} if $C=0$, i.e. if the dual sequence $0 \to {M''}^* \to M^*
  \to {M'}^* \to 0$ is exact. Lemma~\ref{lem6} shows that the
  sequence is dual-exact if and only if $0 \to D(M'') \to D(M) \to D(M') \to
  0$ is exact for a suitable choice of Auslander duals.

  Note also that $C$ is (isomorphic to) a submodule of $\Ext^1(M'',R)$.
\end{Remark}

\subsection{$k$-torsionless modules}
In view of Definition~\ref{def1} and Proposition~\ref{prop5},
we define $k$-torsionless modules as follows:

\begin{Definition} \label{def7}
  A module $M$ is \emph{$k$-torsionless} if $\Ext^i(D(M),R)=0,\,\forall i = 1,
  \dots, k$. Equivalently, $M$ is $k$-torsionless if it is torsionless (for
  $k=1$), if it is reflexive (for $k=2$), resp. if it is reflexive and 
  $\Ext^i(M^*,R)=0,\,\forall i=1,\dots,k-2$ (for $k \geq 3$). By definition,
  \emph{every} module is $0$-torsionless.
\end{Definition}

Note that $M$ is $k$-torsionless if and only if $M_P$ is $k$-torsionless over
$R_P$, $\forall P \in \Spec(R)$ (because both the $\Ext$ functors and the
Auslander dual localize).

\begin{Proposition} \label{prop8}
  Let $0 \to M' \to M \to M'' \to 0$ be an exact sequence, and put
  $C=\Coker(M^*\to {M'}^*)$
  \emph{(cf. Remark after Lemma~\ref{lem6})}.
  \begin{enumerate}
     \item[(a)] If $M'$ and $M''$ are $k$-torsionless and $\grade(C) \geq k$,
                then $M$ is $k$-torsionless.
     \item[(b)] If $M$ is $k$-torsionless, $M''$ is $(k-1)$-torsionless and
                $\grade(C) \geq k-1$, then $M'$ is $k$-torsionless.
     \item[(c)] If $M'$ is $(k+1)$-torsionless, $M$ is $k$-torsionless and
                $\grade(C) \geq k+1$, then $M''$ is $k$-torsionless.
  \end{enumerate}
\end{Proposition}

Recall the definition of $\grade$: $\grade(C) = \inf \{ i \geq 0 \mid 
\Ext^i(C,R) \neq 0 \}$. Since $C$ is finitely generated, we have
$\grade(C) = \depth_I(R) = \inf \{ \depth(R_P) \mid P \in \Supp(C) \}$,
where $I=\Ann_R(C)$.

\begin{proof}
  Using Lemma~\ref{lem6}, we get an exact sequence
  \[
     0 \to C \to D(M'') \to D(M) \to D(M') \to 0,
  \]
  which we can split into short exact sequences
  \[
    0 \to C \to D(M'') \to L \to 0
  \]
  and
  \[
    0 \to L \to D(M) \to D(M') \to 0
  \]
  for suitable $L$. Now the assertions of the Proposition can be proved by
  looking at the corresponding long exact sequences of $\Ext$'s.
\end{proof}

\begin{Remark}
  As $C$ is (isomorphic to) a submodule of $\Ext^1(M'',R)$, the grade condition
  on $C$ holds automatically if it holds for $\Ext^1(M'',R)$. For example
  (see Theorem \ref{thm40} in \S 3), in statements (a) and (b) of the
  Proposition, if $\gdim(M'') \lloc \infty$ then the grade condition on $C$
  follows from the torsionless condition on $M''$. 

  On the other hand, in part (c) the condition on $\grade(C)$ cannot be
  omitted or weakened, as the following example illustrates:
\end{Remark}

\noindent \emph{Example.}
  Let $R$ be any Cohen--Macaulay local ring of dimension $n \geq 1$. Let 
  $\{x_1, \dots, x_n\}$ be a system of parameters. Fix an integer $k,\,0 \leq k
  \leq n-1$. Let $F=R^{k+1},\,x=(x_1, \dots, x_{k+1}) \in F$, and $M=F/Rx$.
  We have an exact sequence:
  \[
     0 \arrow{e} R \arrow{e} R^{k+1} \arrow{e} M \arrow{e} 0.
  \]
  We have $C = \Ext^1(M,R) \isom R/I$, where $I=Rx_1 + \cdots + Rx_{k+1}$.
  Thus $\grade(C)=\depth_I(R)=k+1$; by Proposition \ref{prop8}(c)
  we see that $M$ is $k$-torsionless. (The Proposition applies because
  $R^{k+1}$ and $R$ are free modules.) However, $M$ is \emph{not}
  ($k+1$)-torsionless, or else we would have $\grade(\Ext^1(M,R)) \geq k+2$ by
  Theorem~\ref{thm40} in \S~3 (which applies because $\pd(M)=1$, and in
  particular $\gdim(M) < \infty$).
  \hfill $\Box$

\subsection{Universal pushforward}
Let $M$ be any module (finitely generated with $R$ Noetherian, as always). 
Let $f_1,\dots,f_n \in M^*$ generate $M^*$. Let $\F:M \to R^n$ be the map
$(f_1,\dots,f_n)$. Then $\F$ is injective if and only if $M$ is torsionless;
in fact, $\Ker(\F)=K_M$.

Thus if $M$ is torsionless we get an exact sequence
\begin{equation*} \label{eq:upf}
  0 \arrow{e} M \arrow{e,t}{\F} R^n \arrow{e} N \arrow{e} 0 \tag*{(u.p.f.)}
\end{equation*}
with $N = \Coker(\F)$, and this sequence is also dual-exact ($\F^*$ takes the
canonical basis of $(R^n)^*$ onto $(f_1, \dots,f_n)$; thus $\F^*$ is
surjective). Therefore $\Ext^1(N,R)=0$.

Such an exact and dual-exact sequence, obtained from a system of generators of
$M^*$ -- with $M$ torsionless -- is called a \emph{universal pushforward} of
$M$.

By (c) of Proposition~\ref{prop8}, if $M$ is $k$-torsionless
($k \geq 1$) and \ref{eq:upf} as above is a universal pushforward, then $N$
is $(k-1)$-torsionless. Indeed, $C=\Ext^1(N,R)=0$ in this case.

\vspacec

We conclude this section by showing that $k$-torsionless $\implies$ 
$k^{\text{th}}$ syzygy $\implies$ property $\tilde{S}_k$.

\begin{Definition} \label{def9}
  A module $M$ is a $k^{\text{th}}$ syzygy ($k \geq 1$) if there exists an
  exact sequence 
  \[ 
     0 \to M \to P_0 \to \cdots \to P_{k-1}
  \]
  with $P_j$ projective, $j = 0, \dots, k-1$.

  For $k=0$, \emph{every} module is a $0^{\text{th}}$ syzygy (this is part of
  the definition).
\end{Definition}

\noindent \emph{Example.}
Every dual is a second syzygy. Indeed, let $M$ be any module; then dualize any
projective presentation of $M$ to see that $M^*$ is (at least) a second syzygy.
     
\begin{Definition} \label{def10}
  A module $M$ satisfies property $\tilde{S}_k$ if 
  \[
    \depth(M_P) \geq \min \{k, \depth(R_P) \},\;\;\forall P \in \Spec(R).
  \]
\end{Definition}

$\tilde{S}_k$ is weaker than the better-known Serre condition $S_k$.
A projective module over \emph{any} ring $R$ satisfies $\tilde{S}_k$ for every
$k$; a projective module over $R$ satisfies $S_k$ only if $R$ itself satisfies
$S_k$. From this point of view, property $\tilde{S}_k$ is more like
$k$-torsionlessness and being a $k^{\text{th}}$ syzygy; a projective module is
always $k$-torsionless and a $k^{\text{th}}$ syzygy for every $k$, no matter
what the ring $R$ is.

\begin{Proposition} \label{prop11}
  Let $M$ be any module over a ring $R$. Let $k \geq 0$. Consider the following
  three conditions on $M$: 
    \begin{enumerate}
       \item[(a)] $M$ is $k$-torsionless;
       \item[(b)] $M$ is a $k^{\text{th}}$ syzygy; \emph{and}
       \item[(c)] $M$ satisfies $\tilde{S}_k$.
    \end{enumerate}
  Then \emph{(a) $\implies$ (b) $\implies$ (c)}.
\end{Proposition}

Later in \S 3 we will show that (a), (b) and (c) are equivalent if $M$ has
locally finite Gorenstein dimension, and we will give several other equivalent
conditions.

\begin{proof}
  For $k=0$ all three conditions are automatically true (by definition) for
  every $M$. Let $k \geq 1$.

  (a) $\implies$ (b). $M$ is at least $1$-torsionless, and therefore there is a
  universal push\-forward \ref{eq:upf}. Thus $M$ is at least a first syzygy.
  If $k \geq 2$, then $N$ is $(k-1)$-torsionless; by induction, $N$ is a
  $(k-1)^{\text{st}}$ syzygy, and therefore $M$ is a $k^{\text{th}}$ syzygy.

  (b) $\implies$ (c) \ follows directly from the (otherwise trivial) Depth
  Lemma:

  \begin{DLemma}
    Let $(R,\mm,\kk)$ be a local ring. Let $0 \to K \to L \to M \to 0$ be an
    exact sequence of finitely generated modules. Then
        \[ \depth(K) \geq \min\{\depth(L),\depth(M)+1\}. \]
    If moreover $\depth(L)>\depth(M)$, then $\depth(K)=\depth(M)+1$.
  \end{DLemma}

  The Depth Lemma follows at once from the cohomological characterization of
  depth \cite[Theorem 16.7]{mat}, by looking at the long exact sequence of
  $\Ext^i(\kk,\star)$.
\end{proof}


\section{Gorenstein dimension of modules}

\begin{Definition} \label{def12}
  A module $M$ has Gorenstein dimension $0$
  if $M$ is \emph{reflexive} and
  $\Ext^i(M,R)=\Ext^i(M^*,R)=0,\, \forall i \geq 1$. Equivalently, $\gdim(M)=0$
  if $\Ext^i(M,R) =\Ext^i(D(M),R)=0,\,\forall i \geq 1$.
  \emph{Notation: \mbox{$\gdim(M)=0$}, or $\gdim_R(M)=0$.}
\end{Definition}

\vspaceb

Clearly, if $\gdim(M)=0$, then $\gdim(M^*)=0$ as well, and 
$[\gdim(M) = 0] \iff [\gdim(D(M))=0]$ (recall that $D(D(M)) \proje M$).

\vspaceb

\begin{Lemma} \label{lem13}
  \begin{enumerate}
    \item[(a)] $M$ projective $\implies \gdim(M)=0$.
    \item[(b)] If $0\to K\to L\to M\to 0$ is exact, with $\gdim(M)=0$, then:
               \[
                   [\gdim(K)=0] \iff [\gdim(L)=0].
               \]
    \item[(c)] $\gdim_R(M)=0 \iff\gdim_{R_P}(M_P)=0,\,\forall P\in\Spec(R)$.
  \end{enumerate}
\end{Lemma}

These are all easy to prove.

\begin{Lemma} \label{lem14}
  If $\pd(M) < \infty$ and $\gdim(M)=0$, then $\pd(M)=0$.
\end{Lemma}

\begin{proof}
  If $\pd(M) \leq 1$, we have an exact sequence
  \[
    0 \to P_1 \to P_0 \to M \to 0
  \]
  with $P_j$ projective, $j=0,1$. As $\Ext^1(M,R)=0$, we have
  $0 \to M^* \to P_0^* \to P_1^* \to 0$ exact. Therefore $D(M)\proje 0$, and
  then $M \proje D(D(M)) \proje D(0) \proje 0$, i.e. $M$ is projective.

  If $\pd(M) \geq 2$, consider $0 \to K \to P_0 \to M \to 0$ exact, with $P_0$
  projective; then $\pd(K)<\pd(M)$, and $\gdim(K)=0$ by
  Lemma~\ref{lem13}, (a) and (b). By induction on $\pd(M)$ we see that
  $K$ is projective. But then $\pd(M) \leq 1$, contradicting 
  \mbox{$\pd(M) \geq 2$}.
\end{proof}

\begin{Remark}
  This proof shows that, if $\pd(M) < \infty$ and $\Ext^i(M,R) = 0, \, \forall
  i \geq 1$, then already $\pd(M)=0$. More generally, if $\pd(M) < \infty$ and
  $\Ext^i(M,R)=0, \, \forall i \geq k+1$, then $\pd(M) \leq k$. The analogue
  for Gorenstein dimension is given in Lemma~\ref{lem23}.
\end{Remark}

\begin{Lemma} \label{lem15}
  Let $0 \to M_1 \to M_0 \to M \to 0$ and $0 \to E_1 \to E_0 \to M \to 0$ be
  exact sequences with $M_1$, $M_0$ and $E_0$ having Gorenstein dimension $0$.
  Then $\gdim(E_1)=0$ as well.
\end{Lemma}

\begin{proof}
  Consider the fiber product diagram of $M_0$ and $E_0$ over $M$:
  \begin{equation*}
    \begin{diagram}
      \node[3]{0} \arrow{s} \node{0} \arrow{s}      \\
      \node[3]{E_1} \arrow{e,=} \arrow{s} \node{E_1} \arrow{s}  \\
      \node{0} \arrow{e} \node{M_1} \arrow{e} \arrow{s,=} 
        \node{M_0 \times_M E_0} \arrow{e} \arrow{s} 
        \node{E_0} \arrow{e} \arrow{s} \node{0}     \\
      \node{0} \arrow{e} \node{M_1} \arrow{e} \node{M_0} \arrow{e} \arrow{s} 
        \node{M} \arrow{e} \arrow{s} \node{0}       \\
      \node[3]{0} \node{0} 
    \end{diagram}
  \end{equation*}
   Lemma \ref{lem13}(b) applies twice: first to the top row to give
   $\gdim(M_0 \times_M E_0)=0$, and then to the left column to give
   $\gdim(E_1)=0$.
\end{proof}

\begin{Definition} \label{def16}
  A module $M$ has Gorenstein dimension at most $k$ for some integer $k\geq 0$ 
  \emph{(notation: $\gdim(M) \leq k$)} if there exists an exact sequence 
  \[
    0 \arrow{e} M_k \arrow{e} \cdots \arrow{e} M_0 \arrow{e} M \arrow{e} 0
  \]
  with $\gdim(M_j)=0, \, j=0,\dots,k$.

  If $\gdim(M) \leq k$ for some $k \geq 0$, then we write $\gdim(M) < \infty$.
  Otherwise $\gdim(M) = \infty$.
  If $\gdim(M) < \infty$, we define $\gdim(M)$ as the smallest $k$ such that 
  \mbox{$\gdim(M) \leq k$}.
\end{Definition}

\noindent
Note the perfect similarity with the definition of projective dimension.

\subsection{Characterization of Gorenstein rings in terms of G-dimension}

\begin{Theorem} \label{thm17}
  Assume that $(R,\mm,\kk)$ is local. Let $\dim(R)=n$. Then the following
  conditions are equivalent:
    \begin{enumerate}
      \item[(a)] $R$ is Gorenstein;
      \item[(b)] $\gdim(M) \leq n,\,\forall M$ finitely generated $R$-module;
      \item[(c)] $\gdim(M)<\infty,\,\forall M$ finitely generated $R$-module; 
      \item[(d)] $\gdim(\kk) < \infty$;
      \item[(e)] $\gdim(\kk) = n$.
    \end{enumerate}
\end{Theorem}

\noindent
Compare with the similar characterization of regular local rings in terms of
projective dimension.

\begin{proof}
  We show that (a) $\implies$ (b) $\implies$ (c) $\implies$ (d) $\implies$ (a).
  Then we show that (a)+(b)+(c)+(d) $\implies$ (e); (e) $\implies$ (d) is
  obvious.

  Since (b) $\implies$ (c) $\implies$ (d) are clear, we will show that 
  (a) $\implies$ (b), \mbox{(d) $\implies$ (a)},
  and (a) $+ \cdots +$ (d) $\implies$ (e).

  Recall the various equivalent definitions of Gorenstein local rings
  \cite[Theorem 18.1]{mat}:               \\
  \indent $R$ is Gorenstein $\iff$ $\injdim(R) < \infty$ $\iff$ 
        $\injdim(R)=n$                  \\
  \indent $\iff$ $\Ext^i(\kk,R)=0$ for $i \neq n$ and 
        $\Ext^n(\kk,R) \isom \kk$           \\
  \indent $\iff$ $\Ext^i(\kk,R)=0$ for some $i>n$       \\
  \indent $\iff$ $\Ext^i(\kk,R)=0$ for $i<n$ and $\Ext^n(\kk,R) \isom \kk$.

  \vspaceb

  (a) $\implies$ (b). \ Assume that $R$ is Gorenstein.
      If $n=0$ then $R$ is self-injective, and therefore
      $\Ext^i(M,R)=\Ext^i(D(M),R)=0,\,\forall i \geq 1$ and $\forall M$;
      i.e., all (finitely generated) modules $M$ have Gorenstein dimension $0$.

      Now assume that $n \geq 1$. Let $M$ be a finitely generated $R$-module.
      Let
      \[
        0 \to K \to P_{n-1} \to \cdots \to P_0 \to M \to 0
      \]
      be the beginning of a free resolution of $M$ ($P_j$ free,
      $j=0,\dots,n-1$). Then $\Ext^i(K,R)=\Ext^{i+n}(M,R)=0,\,\forall i\geq 1$,
      for $\injdim(R)=n$.

      Let $F_{\cdot} \to K \to 0$ be a free resolution of $K$.
      As $\Ext^i(K,R)=0,\,\forall i \geq 1$, we have $0 \to K^*\to F_{\cdot}^*$
      exact. In particular, $K^*$ is an $n^{\text{th}}$ syzygy, and then
      $\Ext^i(K^*,R)=0,\,\forall i \geq 1$, for the same reason as for $K$.

      Finally, we show that $K$ is reflexive. As $K$ is at least a first
      syzygy, it is torsionless. Now take $0\to K'\to F\overset{f}{\to} K\to 0$
      exact, with $F$ free. Then $K'$ is an $(n+1)^{\text{st}}$ syzygy; with
      the same proof as for $K$, we get $\Ext^i(K',R) = \Ext^i({K'}^*,R) = 0,\,
      \forall i \geq 1$. In particular, $0 \to K' \to F \to K \to 0$ is
      dual-exact and double-dual-exact. Therefore $f^{**}:F^{**} \to K^{**}$ is
      surjective. As $\sigma_K \circ f = f^{**} \circ \sigma_F$ and $\sigma_F$
      is an isomorphism, we see that $\sigma_K$ is surjective, and therefore an
      isomorphism, as required.

      We have shown that $\gdim(K)=0$, and therefore that $\gdim(M)\leq n$.

  \vspaceb

  (d) $\implies$ (a). \ If $0 \to M_k \to \cdots \to M_0 \to M \to 0$ is an
  exact sequence with $\gdim(M_j)=0,\,j=0,\dots,k$, then $\Ext^{i+k}(M,R)=
  \Ext^i(M_k,R)=0,\,\forall i \geq 1$. In particular, if $\gdim(\kk)<\infty$,
  then $\Ext^i(\kk,R)=0$ for $i > \gdim(\kk)$, so that $R$ is Gorenstein in
  this case.

  \vspaceb

  (a)$+\cdots+$(d) $\implies$ (e). \ By (b), $\gdim(\kk) \leq n$. On the other
  hand, $\Ext^n(\kk,R) \isom \kk \neq 0$, and the proof of (d) $\implies$ (a)
  above shows that $\gdim(\kk) \geq n$.
\end{proof}

\subsection{Gorenstein dimension in exact sequences}

In this subsection we prove the following theorem: 
\begin{Theorem} \label{thm18}
  If $0 \to K \to L \to M \to 0$ is an exact sequence, then:  
  \begin{enumerate}
    \item[(a)] $\gdim(K) \leq \max \{ \gdim(L), \gdim(M)-1 \}$\textup{;}
    \item[(b)] $\gdim(L) \leq \max \{ \gdim(K), \gdim(M) \}$\textup{;}
    \item[(c)] $\gdim(M) \leq 1 + \max \{ \gdim(K), \gdim(L) \}$.
  \end{enumerate}
  
In particular, if two of the three modules $K$, $L$ and $M$ have finite
Gorenstein dimension, then so does the third one.
\end{Theorem}

First we prove several preliminary results. Some of them, in particular 
Theorem~\ref{thm20} and Corollary~\ref{cor22}, are of independent interest.

\begin{Lemma} \label{lem19}
  Let $0 \to K \to L \to M \to 0$ be an exact sequence with $\gdim(M)=0$.
  Then $\gdim(K)=\gdim(L)$.
\end{Lemma}

\begin{proof}
  1. \ Assume first that $\gdim(L)\leq k<\infty$. We show that
  $\gdim(K)\leq k$. The case $k=0$ is covered by Lemma \ref{lem13}(b).
  Assume that $k \geq 1$.

  Let $0 \to F \to T_0 \to L \to 0$ be an exact sequence with $\gdim(T_0)=0$
  and $\gdim(F) \leq k-1$; such an exact sequence exists by the definition of
  $\gdim(L) \leq k$. Consider the commutative diagram 
  \begin{equation*}
    \begin{diagram}
      \node[2]{0} \arrow{s} \node{0} \arrow{s}      \\
      \node[2]{F} \arrow{e,=} \arrow{s} \node{F} \arrow{s}      \\
      \node{0} \arrow{e} \node{S} \arrow{e} \arrow{s} \node{T_0} \arrow{e} 
        \arrow{s} \node{M} \arrow{e} \arrow{s,=} \node{0}       \\
      \node{0} \arrow{e} \node{K} \arrow{e} \arrow{s} \node{L} \arrow{e} 
        \arrow{s} \node{M} \arrow{e} \node{0}               \\
      \node[2]{0}  \node{0}
    \end{diagram}
  \end{equation*}
  with exact rows and columns, where $T_0 \to M$ in the top row is the
  composite of $T_0 \to L$ and $L \to M$, and $S$ is the kernel of $T_0 \to M$.
  (Then the two columns have the same kernel by the Snake Lemma.)

  $\gdim(S)=0$ by Lemma \ref{lem13}(b) applied to the top row; therefore
  $\gdim(K) \leq k$, because $\gdim(F) \leq k-1$.

  \vspaceb

  2. \ Conversely, assume that $\gdim(K) \leq k < \infty$. We show that 
  $\gdim(L) \leq k$, by induction on $k$. The case $k=0$ is covered by 
  Lemma \ref{lem13}(b). Assume that $k \geq 1$.

  Let $0 \to F \to S \to K \to 0$ be an exact sequence with $\gdim(S)=0$ and
  $\gdim(F) \leq k-1$. Let $0 \to D \to P \to M \to 0$ be exact, with $P$
  projective; then $\gdim(D)=0$. Complete the diagram:
    \begin{equation*}
      \begin{diagram}
        \node[2]{0} \arrow{s} \node{0} \arrow{s} \node{0} \arrow{s}  \\
        \node{0} \arrow{e} \node{F} \arrow{e,..} \arrow{s} 
           \node{H} \arrow{e,..} \arrow{s,..} \node{D} \arrow{e} \arrow{s} 
           \node{0}                     \\
        \node{0} \arrow{e} \node{S} \arrow{e} \arrow{s} \node{S \oplus P}
           \arrow{e} \arrow{s,..} \node{P} \arrow{e} \arrow{s} \arrow{sw,..}
           \node{0}                     \\
        \node{0} \arrow{e} \node{K} \arrow{e} \arrow{s} \node{L} \arrow{e} 
           \arrow{s,..} \node{M} \arrow{e} \arrow{s} \node{0}   \\
        \node[2]{0}  \node{0}  \node{0}
      \end{diagram}
    \end{equation*}

  By the inductive hypothesis applied to the top row, $\gdim(H) \leq k-1$;
  therefore $\gdim(L) \leq k$, as required.
\end{proof} 

\begin{Theorem} \label{thm20}
  Let $k \geq 1$. Assume that
  \[
    0 \to K \to M_{k-1} \to \cdots \to M_0 \to N \to 0
  \]
  is exact, with $\gdim(N)\leq k$ and $\gdim(M_j)=0,\,j=0,\dots,k-1$. 

  Then $\gdim(K)=0$.
\end{Theorem}

\noindent
(Compare with the similar statement for projective dimension.)

\begin{proof}
  It suffices to show that [$0 \to K_1 \to M_0 \to N \to 0$ exact, with
  $\gdim(N) \leq k$ and $\gdim(M_0)=0$] $\implies$ [$\gdim(K_1) \leq k-1$] --
  if $k \geq 1$.

  The proof is essentially the same as that of Lemma~\ref{lem15}. Since
  $\gdim(N) \leq k$, there is an exact sequence $0 \to F \to T_0 \to N \to 0$
  with $\gdim(T_0)=0$ and $\gdim(F) \leq k-1$. Now consider the fiber product
  diagram :
  \begin{equation*}
    \begin{diagram}  
      \node[3]{0} \arrow{s} \node{0} \arrow{s}              \\
      \node[3]{F} \arrow{s} \arrow{e,=} \node{F} \arrow{s}      \\
      \node{0} \arrow{e} \node{K_1} \arrow{e} \arrow{s,=}\node{M_0\times_N T_0}
          \arrow{e} \arrow{s} \node{T_0} \arrow{e} \arrow{s} \node{0}   \\
      \node{0} \arrow{e} \node{K_1} \arrow{e} \node{M_0} \arrow{e} \arrow{s} 
          \node{N} \arrow{e} \arrow{s} \node{0}             \\
      \node[3]{0}  \node{0}
    \end{diagram}  
  \end{equation*}
 Lemma \ref{lem19} applies twice: first to the left column, giving 
 $\gdim(M_0 \times_N T_0) \leq k-1$, and then to the top row, giving 
 $\gdim(K_1) \leq k-1$.
\end{proof}

\begin{Corollary} \label{cor21}
  If $\pd(M) < \infty$, then $\gdim(M) = \pd(M)$.
\end{Corollary}

\begin{proof}
  If $\pd(M) < \infty$, we clearly have $\gdim(M) \leq \pd(M)$.

  Now assume that $\gdim(M)=k$ (and $\pd(M) < \infty$). Let
  \[
     0 \to K \to P_{k-1} \to \cdots \to P_0 \to M \to 0
  \]
  be an exact sequence with $P_j$ projective, $j=0, \dots, k-1$.
  Then $\gdim(K)=0$ by Theorem \ref{thm20}. As $\pd(K) < \infty$, 
  Lemma \ref{lem14} shows that $K$ is, in fact, projective -- and therefore
  $\pd(M) \leq k$.

  Alternatively, we could use the Remark after Lemma~\ref{lem14}: we showed in
  the proof of (d) $\implies$ (a) in Theorem~\ref{thm17} that
  $[\gdim(M)=k<\infty] \implies [\Ext^i(M,R)=0 \text{ for all } i>k]$, and this
  together with $\pd(M)<\infty$ implies $\pd(M) \leq k$.
\end{proof}

\begin{Corollary} \label{cor22}
  $\gdim(M) \leq k$ $\iff$ $\gdim_{R_P}(M_P) \leq k,\,\forall P \in\Spec(R)$.

  In other words, $\gdim(M) = \sup \{ \gdim_{R_P}(M_P) \mid P\in\Spec(R) \}$.
\end{Corollary}

\begin{proof}
  $\implies$ is clear, as ``$\gdim(M_j)=0$'' is a local property. Conversely,
  assume that $\gdim_{R_P}(M_P) \leq k,\,\forall P\in\Spec(R)$. Let
  \[
     0 \to K \to M_{k-1} \to \cdots \to M_0 \to M \to 0
  \]
  be exact, with $\gdim(M_j)=0,\,j=0,\dots,k-1$. Then localizing at any $P$ and
  using Theorem \ref{thm20}, we get $\gdim_{R_P}(K_P)=0$;
  thus $\gdim(K)=0$, and therefore $\gdim(M) \leq k$.
\end{proof}

\vspaceb

\noindent
\emph{Proof of Theorem \ref{thm18}.}
  Take a truncated free resolution of $0 \to K \to L \to M \to 0$. 
  For each $k$ we get an exact sequence of $k^{\text{th}}$ syzygies,
  $0 \to K_k \to L_k \to M_k \to 0$. Take, for example, assertion (c). 
  Assume that $\gdim(K) \leq k$ and $\gdim(L) \leq k$. Then
  Theorem~\ref{thm20} gives $\gdim(K_k)=\gdim(L_k)=0$. Thus
  $\gdim(M_k)\leq 1$, and $\gdim(M) \leq 1+k$. 

  (a) and (b) are similar. \hfill       \qed

\subsection{Gorenstein dimension and depth}

In this subsection and the next one, we prove the analogue of the
Auslander--Buchsbaum formula:
\[
  [R \text{ local, $M \neq 0$, and } \gdim(M) < \infty] \implies
  [\gdim(M)+\depth(M)=\depth(R)].
\]

\begin{Lemma} \label{lem23}
  Assume that $\gdim(M) \leq k < \infty$. Then: 
    \begin{enumerate}
           \item[(a)] $\Ext^i(M,R)=0,\,\forall i > k$;
           \item[(b)] If $\Ext^k(M,R)=0$ (and $k \geq 1$), then
                       $\gdim(M) \leq k-1$;
           \item[(c)] If $M \neq 0$, then $\gdim(M) = \max \{ t \geq 0 \mid
                       \Ext^t(M,R) \neq 0 \}$.
      \end{enumerate}
\end{Lemma}

\begin{Remark}
  Say $M \neq 0$ and $\gdim(M) < \infty$. Compare: $\gdim(M) = \max \{ t \geq 0
  \mid \Ext^t(M,R) \neq 0 \}; \; \grade(M) = \min \{ t \geq 0 \mid 
  \Ext^t(M,R) \neq 0 \}$. In particular, $\grade(M) \leq \gdim(M)$.
\end{Remark}

\noindent
\emph{Proof of Lemma \ref{lem23}. }
 \ We have already seen (a) in the proof of (d) $\implies$ (a) in
Theorem~\ref{thm17}, and (c) clearly follows from (a) and (b). 

We prove (b). First assume that $\gdim(M) \leq 1$ and $\Ext^1(M,R)=0$; we show
that $\gdim(M)=0$. Let $0 \to M_1 \to M_0 \to M \to 0$ be an exact sequence
with $\gdim(M_j)=0,\,j=0,1$. As $\Ext^1(M,R)=0$, the sequence is dual-exact,
and therefore we have an exact sequence $0 \to D(M) \to D(M_0) \to D(M_1) \to
0$ for suitable choices of Auslander duals (Remark after Lemma~\ref{lem6}).
$\gdim(M_j)=0 \implies \gdim(D(M_j))=0,\,j=0,1$; therefore $\gdim(D(M))=0$, and
finally $\gdim(M)=0$.

Now let $\gdim(M) \leq k$ with $k \geq 2$, and assume that $\Ext^k(M,R)=0$.
Let $0 \to K \to M_{k-2} \to \cdots \to M_0 \to M \to 0$ be an exact sequence
with $\gdim(M_j)=0,\,j=0,\dots,k-2$. Then $\gdim(K) \leq 1$ and $\Ext^1(K,R) =
\Ext^k(M,R) = 0$. By the case $k=1$, discussed above, we get $\gdim(K)=0$, and
therefore $\gdim(M) \leq k-1$.  \hfill      \qed

\begin{Corollary} \label{cor24}
  If $\gdim(M)=0$ and $0 \to M \to R^n \to N \to 0$ is a universal pushforward,
  then $\gdim(N)=0$.
\end{Corollary}

\begin{proof}
  $\Ext^1(N,R)=0$, and on the other hand $\gdim(N) \leq 1$; part (b) of the
  Lemma gives $\gdim(N)=0$.
\end{proof}

\begin{Lemma} \label{lem25}
  (a) Let $0 \to M_1 \to M_0 \to M \to 0$ be an exact sequence with $M_1$ and
  $M_0$ reflexive and $\Ext^1(M_0,R)=0$. Then $\Ext^1(M,R)^*=0$.

  (b) Let $(R,\mm,\kk)$ be local, with $\depth(R)=0$. If $F$ is a finitely
  generated $R$-module with $F^*=0$, then $F=0$.
\end{Lemma}

\begin{Corollary} \label{cor26}
  If $\gdim(M) < \infty$ and $(R,\mm,\kk)$ is local with $\depth(R)=0$, then
  $\gdim(M)=0$.
\end{Corollary}

\begin{proof}
  It suffices to consider the case $\gdim(M) \leq 1$. Let $0 \to M_1 \to M_0
  \to M \to 0$ be exact, with $\gdim(M_j)=0,\,j=0,1$. Then part (a) of 
  Lemma~\ref{lem25} gives $\Ext^1(M,R)^*=0$, and part (b) gives
  $\Ext^1(M,R)=0$. Therefore $\gdim(M)=0$ by Lemma \ref{lem23}(b).
\end{proof}

\noindent
\emph{Proof of Lemma~\ref{lem25}.}
\ (a) \ Dualizing the given sequence we get $0 \to M^* \to M_0^* \to M_1^* \to
\Ext^1(M,R) \to 0$ exact. Dualizing again, we get $0 \to \Ext^1(M,R)^* \to
M_1^{**} \to M_0^{**}$ exact. But $M_1^{**} \to M_0^{**}$ is the same as the
injective map $M_1 \to M_0$ (up to canonical isomorphisms $\sigma_{M_1}$ and
$\sigma_{M_0}$), and therefore $\Ext^1(M,R)^*=0$.

\ (b) \ If $F \neq 0$, then $F/\mm F$ is a non-zero $\kk$-vector space.
Let $p:F \to F/\mm F$ be the canonical surjection, and let $\pi : F/\mm F \to
\kk$ be any non-zero (and therefore surjective) $\kk$-linear map. Finally,
$\depth(R)=0 \implies \exists\, \alpha : \kk \to R$, a non-zero (and therefore
injective) $R$-linear map. Then $u = \alpha \circ \pi \circ p \in F^*$ is
non-zero, and therefore $F^* \neq 0$.  \hfill       \qed

\begin{Lemma} \label{lem27}
  Let $(R,\mm,\kk)$ be local, with $\depth(R)=0$. Let $E$ be any nonzero,
  finitely generated $R$-module. Then $\depth(E^*)=0$.

  In particular, if $M \neq 0$ is any \emph{reflexive} module (for example, if 
  $\gdim(M)=0$), then $\depth(M)=0$.
\end{Lemma}

\begin{proof} 
  As in the proof of Lemma \ref{lem25}(b), we can find a surjective
  homomorphism $E \to \kk$. Dualizing we get an injective homomorphism 
  $\kk^* \to E^*$. $\kk^*$ is nonzero (because $\depth(R)=0$) and of finite
  length, and therefore $\depth(E^*)=0$.
\end{proof}

\begin{Lemma} \label{lem28}
  Let $(R,\mm,\kk)$ be local, of arbitrary depth, and let $\gdim(M)=0$ with
  $M \neq 0$. Then $\depth(M)=\depth(R)$.
\end{Lemma}

\begin{proof} 
  Put $\depth(R)=d$.
  
  $M$ is $k$-torsionless for every $k \geq 0$. In particular, $M$ is
  $d$-torsionless. By Proposition~\ref{prop11}, $M$ satisfies $\tilde{S}_d$,
  and in particular $\depth(M) \geq d$.

  Now we show that $\depth(M)=d$. The case $d=0$ is covered by Lemma
  \ref{lem27}; therefore we may assume that $d \geq 1$. We will show
  that $\Ext^d(\kk,M) \neq 0$ (which means that $\depth(M) \leq d$).

  We have $\gdim(M^*)=0$. Consider the beginning of a minimal free resolution
  of $M^*$:
  \[
     F_1 \arrow{e,t}{u} F_0 \arrow{e,t}{f} M^* \arrow{e} 0.
  \]

  Put $K_1 = \Ker(f),\; K_2 = \Ker(u)$; thus we have exact sequences:
  \[
     0 \to K_2 \to F_1 \to K_1 \to 0  \quad \text{and} \quad
     0 \to K_1 \to F_0 \to M^* \to 0,
  \]
  and $\gdim(K_j)=0,\,j=1,2$. Dualizing, we get exact sequences
  \[
     0 \to M \to F_0^* \to K_1^* \to 0 \quad \text{and} \quad
     0 \to K_1^* \to F_1^* \to K_2^* \to 0
  \]
  (using $M^{**} \isom M$).

  $\gdim(K_1^*)=0 \implies \depth(K_1^*) \geq d \implies 
  \Ext^{d-1}(\kk,K_1^*)=0$; therefore we have an exact sequence
  \[
     0 \arrow{e} \Ext^d(\kk,M) \arrow{e} \Ext^d(\kk,F_0^*)
       \arrow{e,t}{\alpha} \Ext^d(\kk,K_1^*).
  \]

  We show that $\alpha=0$; then $\Ext^d(\kk,M)=\Ext^d(\kk,F_0^*) \neq 0$
  (because $\depth(F_0^*)=\depth(R)=d$), as required.

  To show that $\alpha=0$, consider the commutative diagram
    \begin{equation*}
    \begin{diagram}          
       \node{\Ext^d(\kk,F_0^*)} \arrow[2]{e,t}{\gamma} \arrow{se,b}{\alpha}
        \node[2]{\Ext^d(\kk,F_1^*)} \\
       \node[2]{\Ext^d(\kk,K_1^*)} \arrow{ne,b}{\beta}
    \end{diagram}          
      \end{equation*}
  $\Ker(\beta) = \Ext^{d-1}(\kk,K_2^*)=0$, so that $\beta$ is injective. On the
  other hand, $\gamma$ is induced by $u^* : F_0^* \to F_1^*$, and therefore
  $\gamma=0$. Consequently $\alpha=0$, as advertised.
\end{proof}

\begin{Theorem}[The Auslander--Bridger Formula]  \label{ABF}
  Let $(R,\mm,\kk)$ be a local Noetherian ring, and let $M \neq 0$ be a
  finitely generated $R$-module. If $\gdim(M) < \infty$, then
  \[  \gdim(M) + \depth(M) = \depth(R).  \]
\end{Theorem}

\begin{proof} 
  Induction on $k=\gdim(M)$. If $k=0$, the result is contained in Lemma
  \ref{lem28}.

  If $k \geq 1$, then $\depth(R) \geq 1$ by Corollary \ref{cor26}. Put
  $\depth(R)=d$.

  Assume first that $k=1$. Then we have an exact sequence
  \[
     0 \to K \to F \to M \to 0
  \]
  with $F$ free and $\gdim(K)=0$. Then $\depth(K)=\depth(F)=d$. By the Depth
  Lemma, $\depth(M) \geq d-1$. The crux of the whole proof is to show that 
  $\depth(M)=d-1$ in this case. We will prove this later, in the next
  subsection, after we will have studied the behavior of Gorenstein dimension
  under reduction modulo a regular element.

  If $k \geq 2$, we have an exact sequence 
  \[
     0 \to K \to F \to M \to 0
  \]
  with $F$ free and $\gdim(K)=k-1$. By induction, $\depth(K)=d-(k-1)$, and then
  $\depth(M) = \depth(K)-1 = d-k$ by the Depth Lemma, for $\depth(F) = d >
  \depth(K)$ (so that $\depth(K) \geq \depth(M)+1$, and in particular
  $\depth(F) > \depth(M)$).
\end{proof}

\begin{Remarks}
  (1) \ The proof in \cite[pp. 118--119]{ab} contains a serious mistake.
  Namely, when $\gdim(M) \geq 1$, the authors consider an exact sequence
  \[
      0 \to K \to P \to M \to 0
  \]
  with $P$ projective, and then claim that ``It is well-known and easily
  proved that $\depth(K)=\depth(M)+1$.'' Unfortunately, this formula is not
  necessarily true when $\depth(K)=\depth(P)$; this is exactly why the case
  $\gdim(M)=1$ is the most difficult. 

  At that point in the memoir, Auslander and Bridger had already included
  reduction modulo a regular element in the theory. They used it for the case
  $\gdim(M)=0$. As we have seen, that case is, in fact, elementary (i.e. does
  not require reduction modulo a regular element). I was unable to find a
  similar elementary proof for the case $\gdim(M)=1$. Compare with the much
  easier proof of the usual Auslander--Buchsbaum formula for projective
  dimension in \cite[p.~155]{mat}.

  (2) \ Let $R$ be a ``global'' (i.e. not necessarily local) Noetherian ring.
  Assume that $\gdim(M) \lloc \infty$, i.e. $\gdim_{R_P}(M_P) < \infty,\,
  \forall P \in \Spec(R)$. If $\dim(R) < \infty$ ($\dim(R)$ is the Krull
  dimension of $R$), then $\gdim(M) \leq \dim(R)$; this follows from the
  Auslander--Bridger formula. However, if $\dim(R)=\infty$ it is possible to
  have $\gdim(M)\lloc\infty$ but $\gdim(M)=\infty$. For example, let $\kk$ be
  a field, and let $A=\kk[X_1,\dots,X_n,\dots]$ be the polynomial ring over
  $\kk$ in countably many variables. Let $1 = m_1 < m_2 < \cdots$ be positive
  integers such that $m_2-m_1 < m_3-m_2 < \cdots$. Let $\mm_i$ be the prime
  ideal of $A$ generated by $\{X_j \mid m_i \leq j < m_{i+1} \}$, let $S$ be
  the complement of $\union_i \mm_i$ in $A$, and let $R=S^{-1}A$. Then $R$ is a
  Noetherian ring of infinite Krull dimension \cite[p. 203]{nagata}. In this
  example, $\kk$ is naturally a finitely generated $R$-module, and $\gdim(\kk)
  \lloc \infty$, but $\gdim(\kk) = \infty$.
\end{Remarks}

\begin{Corollary} \label{cor30}
  If $\gdim(M) \lloc \infty$, then 
  \[ \grade(\Ext^i(M,R)) \geq i,\, \forall i \geq 1.  \]
\end{Corollary}

\begin{proof} \ \ 
  $\grade(\Ext^i(M,R)) = \min \{\depth(R_P) \mid P \in\Supp(\Ext^i(M,R))\}$.\\
  If $\depth(R_P) < i$, then $\gdim_{R_P}(M_P) \leq \depth(R_P) < i$ by the
  Auslander--Bridger formula, so that $\Ext_{R_P}^i(M_P,R_P)=0$ by
  Lemma \ref{lem23}(a). Therefore $\Ext^i(M,R)_P=0$, i.e. 
  $P \not \in \Supp(\Ext^i(M,R))$.
\end{proof}

\subsection{Regular elements and Gorenstein dimension}

In this subsection we complete the proof of the Auslander--Bridger formula.
Before we can do that, however, we must bring regular elements into the
picture. In particular, we need to study the behavior of Gorenstein dimension
under reduction modulo a regular element.

\begin{Notation}
  Throughout this subsection, fix a ring $R$ and a non-unit $x \in R$. Write
  $\bar{R}$ for $R/xR$, $\bar{M}$ for $M/xM$, etc. Note that $\bar{M}^* =
  \Hom_{\bar{R}}(\bar{M},\bar{R})$, while $\overline{M^*} = \Hom_R(M,R) \tens_R
  \bar{R}$. Similarly with $D_{\bar{R}}(\bar{M})$ and $\overline{D_R(M)}$.
\end{Notation}

\vspaceb

Consider the functors $\star \, \tens_R \bar{R}$ and
$\Hom_{\bar{R}}(\star,\bar{R})$ and their composite,
$\Hom_R(\star,\bar{R})$; note that we have an obvious functorial isomorphism 
  \[  \Hom_{\bar{R}}(\bar{M},\bar{R}) \isom \Hom_R(M,\bar{R}). \]

In this case the abstract Leray spectral sequence becomes
\[
   \Ext_{\bar{R}}^i(\Tor_j^R(M,\bar{R}),\bar{R}) \, \Rightarrow \,
        \Ext_R^{i+j}(M,\bar{R}).
\]
In particular, if $x$ is both $R$-regular and $M$-regular, then
$\Tor_j^R(M,\bar{R})=0,\,\forall j \geq 1$, and we get canonical isomorphisms
\[
   \Ext_{\bar{R}}^i(\bar{M},\bar{R}) \isom \Ext_R^i(M,\bar{R}),\,
            \forall i \geq 0.
\]

Of course, these isomorphisms can be proved easily (when $x$ is $R$-regular and
$M$-regular) without using spectral sequences. Here is a quick argument for the
benefit of the graduate student:

First, $\Hom_{\bar{R}}(\bar{M},\bar{R}) \isom \Hom_R(M,\bar{R})$ is clear.

Next, let $0 \to K \to F \to M \to 0$ be an exact sequence with $F$ free of
finite rank. Then $x$ is also $K$-regular, and the sequence $0 \to \bar{K} \to
\bar{F} \to \bar{M} \to 0$ is exact (by the Snake Lemma). We have a commutative
diagram with exact rows:
\begin{equation*}
  \begin{diagram}
    \node{\Hom_R(F,\bar{R})} \arrow{e} \arrow{s,r}{\vertisom}
        \node{\Hom_R(K,\bar{R})} \arrow{e} \arrow{s,r}{\vertisom}
        \node{\Ext^1_R(M,\bar{R})} \arrow{e} \arrow{s,..} \node{0} \\
    \node{\Hom_{\bar{R}}(\bar{F},\bar{R})} \arrow{e}
    \node{\Hom_{\bar{R}}(\bar{K},\bar{R})} \arrow{e}
    \node{\Ext^1_{\bar{R}}(\bar{M},\bar{R})} \arrow{e} \node{0}
  \end{diagram}
\end{equation*}
and $\Ext^1_{\bar{R}}(\bar{M},\bar{R}) \isom \Ext^1_R(M,\bar{R})$ follows.

Finally, for $i \geq 2$ we have by induction on $i$:
\[
  \Ext^i_{\bar{R}}(\bar{M},\bar{R}) \isom \Ext^{i-1}_{\bar{R}}(\bar{K},\bar{R})
  \isom \Ext^{i-1}_R(K,\bar{R}) \isom \Ext^i_R(M,\bar{R}).
\]

\vspacec

Since the Auslander dual commutes with base change (Remark \#5 after
Proposition~\ref{prop4}), we have:
\begin{Proposition} \label{prop31}
  If $x$ is any element of $R$ and $M$ is any $R$-module, then
  \[  D_{\bar{R}}(\bar{M}) \proje \overline{D_R(M)}. \]
\end{Proposition}
\hfill \qed

\begin{Corollary} \label{cor32}
  If $x$ is $R$-regular and $\gdim(M)=0$, then $x$ is $M$-regular and
  $\gdim_{\bar{R}}(\bar{M})=0$.
\end{Corollary}

\begin{proof} 
  First we show that $x$ is $M$-regular. Indeed, $M$ is $1$-torsionless, and
  therefore a first syzygy (Proposition~\ref{prop11}). Since $x$ is $R$-regular
  and $M$ is isomorphic to a submodule of a free module, $x$ is $M$-regular.

  Therefore $\Ext_{\bar{R}}^i(\bar{M},\bar{R}) \isom \Ext_R^i(M,\bar{R}),\,
  \forall i \geq 0$. From $0 \to R \overset{\cdot x}{\to} R \to \bar{R} \to 0$
  exact we get
  \[
    \Ext_R^i(M,R) \arrow{e} \Ext_R^i(M,\bar{R}) \arrow{e} 
        \Ext_R^{i+1}(M,R)
  \]
  exact. Since $\gdim(M)=0$, we get $\Ext_{\bar{R}}^i(\bar{M},\bar{R}) =
  \Ext_R^i(M,\bar{R}) = 0,\, \forall i \geq 1$.

  On the other hand, $\gdim_R(M)=0 \implies \gdim_R(D(M))=0$; by what we have
  just proved, $\Ext_{\bar{R}}^i(\overline{D_R(M)},\bar{R}) = 0,\, \forall i
  \geq 1$. But $\overline{D_R(M)} \proje D_{\bar{R}}(\bar{M})$. Thus
  \[
     \Ext_{\bar{R}}^i(\bar{M},\bar{R}) =
     \Ext_{\bar{R}}^i(D_{\bar{R}}(\bar{M}),\bar{R}) = 0,\, \forall i \geq 1
  \]
  -- that is, $\gdim_{\bar{R}}(\bar{M}) = 0$.
\end{proof}

\begin{Corollary} \label{cor33}
  If $x$ is $R$-regular and $M$-regular and $\gdim(M) < \infty$, then \\
  $\gdim_{\bar{R}}(\bar{M}) \leq \gdim_R(M)$.
  If moreover $x \in J(R)$ (the Jacobson radical of $R$), then
  $\gdim_{\bar{R}}(\bar{M}) = \gdim_R(M)$.
\end{Corollary}

\begin{proof} 
  Induction on $k=\gdim_R(M)$; the case $k=0$ is covered by the previous
  Corollary. Assume that $k \geq 1$.

  Consider an exact sequence 
  \[  0 \to K \to F \to M \to 0  \]
  with $F$ free and $\gdim(K)=k-1$. Since $x$ is $M$-regular and $F$-regular,
  it is also $K$-regular, and we have an exact sequence of $\bar{R}$-modules 
  \[  0 \to \bar{K} \to \bar{F} \to \bar{M} \to 0.  \]

  By induction, $\gdim_{\bar{R}}(\bar{K}) \leq k-1$, and therefore
  $\gdim_{\bar{R}}(\bar{M}) \leq k$.

  To show equality when $x \in J(R)$, it suffices to show that
  $\Ext_{\bar{R}}^k(\bar{M},\bar{R}) \neq 0$ in that case (see
  Lemma \ref{lem23}(a)).

  Equivalently, we will show that $\Ext_R^k(M,\bar{R}) \neq 0$. We have an
  exact sequence:
  \[
    \Ext_R^k(M,R) \arrow{e,t}{\cdot x} \Ext_R^k(M,R) \arrow{e}
       \Ext_R^k(M,\bar{R}).
  \]
  As $x \in J(R)$, Nakayama's Lemma shows that $[ \Ext_R^k(M,\bar{R})=0 ]
  \implies [ \Ext_R^k(M,R)=0 ]$; but if that happens, then $\gdim_R(M)<k$,
  contradiction.
\end{proof}

\vspacec

Now we are ready to finish the proof of the Auslander--Bridger formula.

\vspaceb

\noindent
\emph{Proof of Theorem \ref{ABF}: Conclusion.}

\vspaceb

The only thing left to prove is: If $(R,\mm)$ is local, $d=\depth(R)$, and
$\gdim_R(M)=1$, then $\depth(M)=d-1$.

We prove this by induction on $d$. The case $d=0$ is vacuously true (in that
case there are \emph{no} modules $M$ with $\gdim(M)=1$). Let $d \geq 1$.

We have already proved that $\depth(M) \geq d-1$. By way of contradiction,
assume that $\depth(M) \geq d$. Then there exists an element $x \in \mm$ which
is both $R$-regular and $M$-regular. We have: $\depth(\bar{R})=d-1$ and
$\gdim_{\bar{R}}(\bar{M}) = \gdim_R(M) = 1$. By induction,
$\depth_{\bar{R}}(\bar{M}) = d-2$; equivalently, $\depth_R(\bar{M}) = d-2$.

But now $\depth(M)=d$ and $\depth(\bar{M})=d-2$ taken together contradict the
Depth Lemma applied to the exact sequence 
$0 \to M \overset{\cdot x}{\to} M \to \bar{M} \to 0$.
\hfill \qed

\vspaceb

\noindent
{\bf Exercise.} \ (cf. \cite[Theorem 16.9]{mat} for the case of projective
dimension.)

Let $R$ be a Noetherian ring, and let $M,N$ be finitely generated $R$-modules
with $M \neq 0$, $\grade(M)=k>0$, and $\gdim(N)=l<k$.

Show that $\Ext^i(M,N) = 0, \, \forall i < k-l$.

\vspacea

I include this exercise to further illustrate the theme of results which remain
true if one replaces projective dimension with Gorenstein dimension (sometimes
at the price of more delicate proofs). 

\vspacea

\emph{Hint:} Induction on $l=\gdim(N)$. The case $l=0$ is the most serious.
(The case with $\pd(N)=0$ is trivial by comparison.)

Assume that $\gdim(N)=0$. First show that $\Hom(M,N)=0$: as $\grade(M)=k>0$,
$\exists\, x \in \Ann(M)$ with $x$ $R$-regular. Then $x$ is also $N$-regular,
and $\Hom(M,N)=0$ follows easily. Then take a universal pushforward
$0 \to N \to F \to N' \to 0$ (with $F$ free and $\gdim(N')=0$; see Corollary
\ref{cor24}). As $\Ext^i(M,F)=0, \, \forall i<k$, we get 
$\Ext^i(M,N)=0, \, \forall i<k$, by bootstrapping.


\section{$k^{\text{th}}$ syzygies of finite Gorenstein dimension}

In this section we study the relationship between $k$-torsionless, being a
$k^{\text{th}}$ syzygy and condition $\tilde{S}_k$ for modules of locally
finite Gorenstein dimension. Proposition~\ref{prop11} in \S 1 gives
implications which hold for an arbitrary module $M$. (We will refine that
result in Proposition~\ref{prop36} below.)

\vspaceb

First we introduce $k$-torsionfreeness:

\begin{Definition} \label{def34}
  Let $(R,\mm)$ be a local ring, and let $M$ be an $R$-module. Fix an integer
  $k \geq 0$. $M$ is \emph{$k$-torsionfree} if every $R$-regular sequence of
  length at most $k$ is also $M$-regular. \emph{(Note that ``1-torsionfree''
  is usually called ``torsionfree''.)}
\end{Definition}

\begin{Definition} \label{def35}
  Let $R$ be any ring (not necessarily local). A module $M$ is locally
  $k$-torsionfree if $M_P$ is $k$-torsionfree over $R_P$, $\forall P \in
  \Spec(R)$.
\end{Definition}

As usual, every module is (locally) $0$-torsionfree.

\begin{Remark}
  Auslander and Bridger use the term ``$k$-torsionfree'' for what I call
  ``$k$-torsionless'' (and they do not use any name for what I call
  ``$k$-torsionfree''). It seems to me that the meaning of ``torsionfree''
  as given in Definition~\ref{def34} above is standard in algebra -- and
  therefore I use ``torsionless'' for what \cite{ab} calls ``torsionfree''.
\end{Remark}

\vspaceb

Local $k$-torsionfreeness lies between being a $k^{\text{th}}$ syzygy and
property $\tilde{S}_k$. In fact, we have the following refinement of
Proposition~\ref{prop11}:

\begin{Proposition} \label{prop36}
  Let $M$ be any module over a ring $R$. Let $k \geq 0$. Consider the following
  conditions on $M$:
    \begin{enumerate}
       \item[(a)] $M$ is $k$-torsionless;
       \item[(b)] $M$ is a $k^{\text{th}}$ syzygy;
       \item[(c)] $\exists \: 0 \to M \to M_0 \to \dots \to M_{k-1}$ exact,
                         with $\gdim(M_j)=0, \, j=0,\dots,k-1$;
       \item[(d)] $M$ is locally $k$-torsionfree;
       \item[(e)] $M$ satisfies $\tilde{S}_k$.
    \end{enumerate}
  Then \emph{(a) $\implies$ (b) $\implies$ (c) $\implies$ (d)
  $\implies$ (e)}.
\end{Proposition}

\emph{Proof.}
  We have seen that (a) $\implies$ (b) in Proposition \ref{prop11}, and (b)
  $\implies$ (c) and (d) $\implies$ (e) are clear.

  We prove (c) $\implies$ (d). By localizing at any prime $P$, it suffices to
  show that (c) $\implies$ [$M$ is $k$-torsionfree] over a \emph{local} ring
  $R$. We do induction on $k$, the case $k=0$ being trivial. If $k \geq 1$, let
  $N=\Coker(M \to M_0)$; then we have a short exact sequence
  \[ 0 \to M \to M_0 \to N \to 0  \]
  where $N$ satisfies (c) for $k-1$, and therefore $N$ is $(k-1)$-torsionfree
  by induction. $M_0$ is $l$-torsionfree for every $l \geq 0$; this follows
  easily from Corollary~\ref{cor32}. Therefore the conclusion follows from the
  following Lemma:

\begin{Lemma} \label{lem37}
  If $0 \to M \to B \to N \to 0$ is a short exact sequence with $B$
  $k$-torsionfree and $N$ $(k-1)$-torsionfree, then $M$ is $k$-torsionfree.

  \emph{Compare this with Proposition \ref{prop8}(b).}
\end{Lemma}

\begin{proof}
  If $k \geq 1$ and $x_1$ is $R$-regular, then $x_1$ is $B$-regular and
  therefore $M$-regular as well. Thus $M$ is at least $1$-torsionfree in this
  case.

  If $k \geq 2$ and $x_1,\dots,x_s$ is an $R$-regular sequence with $2 \leq s
  \leq k$, then $x_1$ is $N$-regular, $B$-regular and $M$-regular, and we have
  an exact sequence of $\bar{R}$-modules,
  \[ 0 \to \bar{M} \to \bar{B} \to \bar{N} \to 0,  \]
  where $\bar{R}=R/xR$, $\bar{M}=M/xM$, etc.

  As $B$ is $k$-torsionfree over $R$ and $x_1$ is $R$-regular, we see that 
  $\bar{B}$ is $(k-1)$-torsionfree over $\bar{R}$. Similarly, $\bar{N}$ is
  $(k-2)$-torsionfree over $\bar{R}$. By induction, $\bar{M}$ is
  $(k-1)$-torsionfree over $\bar{R}$, and in particular $x_2,\dots,x_s$ is
  $\bar{M}$-regular. But then $x_1,x_2,\dots,x_s$ is $M$-regular, as required.
\end{proof}

\vspaceb

When $\gdim(M) \lloc \infty$, all the conditions in Proposition~\ref{prop36}
are equivalent. The cycle of implications is closed via yet another equivalent
condition, the prototype of which we have already seen in
Corollary~\ref{cor30}.

\begin{Proposition} \label{prop38}
  Fix an integer $k \geq 0$. Assume that $\gdim(M) \lloc \infty$ and that $M$
  satisfies $\tilde{S}_k$. 

  Then $\grade(\Ext^i(M,R)) \geq i+k, \; \forall i \geq 1$.
\end{Proposition}

\noindent The proof is essentially the same as that of Corollary \ref{cor30}:

\begin{proof}
  Fix $i \geq 1$, and take $P \in \Spec(R)$ with $\depth(R_P) < i+k$. Then we
  must show that $\Ext_{R_P}^i(M_P,R_P)=0$, i.e. that $P \not\in
  \Supp(\Ext^i(M,R))$. We have $\gdim_{R_P}(M_P) = \depth(R_P) - \depth(M_P)
  \leq \depth(R_P) - \min \{ k;\depth(R_P) \} = \max \{0; \depth(R_P)-k \}$.
  Since $\depth(R_P) < i+k$ and $i \geq 1$, we have in any case
  $\gdim_{R_P}(M_P) < i$, and therefore $\Ext_{R_P}^i(M_P,R_P) = 0$ by
  Lemma \ref{lem23}(a).
\end{proof}

\begin{Proposition} \label{prop39}
  Fix an integer $k \geq 0$. Assume that $\gdim(M)=g<\infty$, and that
  $\grade(\Ext^i(M,R)) \geq i+k$ for $1 \leq i \leq g$. Then $M$ is
  $k$-torsionless.
\end{Proposition}

\begin{proof}
  Induction on $g \geq 0$. For $g=0$ the hypothesis is just that 
  $\gdim(M)=0$ -- but in that case $M$ {\bf is} $k$-torsionless for every 
  $k \geq 0$. 

  Assume that $g \geq 1$. Let 
  \[ 0 \to M' \to P \to M \to 0 \]
  be an exact sequence with $P$ projective and $\gdim(M')=g-1$. For
  $1 \leq i \leq g-1$, we have
  \[ \grade(\Ext^i(M',R)) = \grade(\Ext^{i+1}(M,R)) \geq i+k+1; \]
  by the inductive hypothesis, $M'$ is $(k+1)$-torsionless. $P$ is
  projective, and in particular $(k+1)$-torsionless; therefore $M$ is
  $k$-torsionless by Proposition \ref{prop8}(c). Note that the hypothesis 
  $\grade(\Ext^1(M,R)) \geq 1+k$ is used here, while the other grade conditions
  are used to show (by induction) that $M'$ is $(k+1)$-torsionless.
\end{proof}

\vspaceb

Putting all these results together, we get:

\begin{Theorem} \label{thm40}
  Fix $k \geq 0$, and assume that $\gdim(M) \lloc \infty$ (where $M$ is a 
  module over any ring $R$, not necessarily local).

  Then the following conditions on $M$ are equivalent:
    \begin{enumerate}
       \item[(a)] $M$ is $k$-torsionless;
       \item[(b)] $M$ is a $k^{\text{th}}$ syzygy;
       \item[(c)] $\exists \, 0 \to M \to M_0 \to \dots \to M_{k-1}$
                  exact, with $\gdim(M_j)=0, \, j=0,\dots,k-1$;
       \item[(d)] $M$ is locally $k$-torsionfree;
       \item[(e)] $M$ satisfies $\tilde{S}_k$;
       \item[(f)] $\grade(\Ext^i(M,R)) \geq i+k, \, \forall i \geq 1$.
    \end{enumerate}
\end{Theorem}

In particular, all these conditions are equivalent if $M$ is any module over
a Gorenstein ring $R$, or if $\pd(M) \lloc \infty$.  \hfill \qed

\vspacec

The equivalence of (a) and (b) holds in more general conditions. A typical
situation where the result below can be used is that of modules which have 
(locally) finite projective (or Gorenstein) dimension on the punctured
spectrum of a local ring of sufficiently large depth.

\begin{Theorem} \label{thm41}
  Fix $k \geq 0$. Let $R$ be a ring and $M$ an $R$-module. Assume that 
  \[
    [P \in \Spec(R) \text{ and } \depth(R_P) \leq k-2] \implies 
    [\gdim_{R_P}(M_P) < \infty].
  \]

  Then the following conditions on $M$ are equivalent:
    \begin{enumerate}
       \item[(a)] $M$ is $k$-torsionless;
       \item[(b)] $M$ is a $k^{\text{th}}$ syzygy.
    \end{enumerate}
\end{Theorem}

\begin{proof}
  We already know that (a) $\implies$ (b), with no conditions on G-dimension.
  Conversely, assume that $M$ is a $k^{\text{th}}$ syzygy and that 
  $\gdim_{R_P}(M_P) < \infty$ whenever $\depth(R_P) \leq k-2$. We show that $M$
  is $k$-torsionless by induction on $k$. For $k=0$ there is nothing to prove,
  and the case $k=1$ is easy. Now let $k \geq 2$. Since $M$ is a
  $k^{\text{th}}$ syzygy, there is an exact sequence 
    \[   0 \to M \to F \to N \to 0   \]
  with $F$ free and $N$ a $(k-1)^{\text{st}}$ syzygy. By induction, $N$ is
  $(k-1)$-torsionless; if we can show that $\grade(\Ext^1(N,R)) \geq k-1$, then
  Proposition \ref{prop8}(a) gives that $M$ is $k$-torsionless, as required.

  To finish the proof, we show that $\grade(\Ext^1(N,R)) \geq k-1$.
  Let $P \in \Spec(R)$ with $\depth(R_P) < k-1$. We must show that $P \notin
  \Supp(\Ext^1(N,R))$. But this is clear, since
  \begin{multline*}
    \depth(R_P) < k-1 \implies \gdim_{R_P}(M_P) < \infty  \\
    \implies \gdim_{R_P}(N_P) < \infty 
      \implies \grade(\Ext_{R_P}^1(N_P,R_P)) \geq 1+(k-1),
  \end{multline*}
  by Proposition~\ref{prop38} ($N_P$ is a $(k-1)^{\text{st}}$ syzygy, and
  therefore satisfies $\tilde{S}_{k-1}$); then, in fact, $\Ext^1_{R_P}(N_P,R_P)
  =0$, because the grade of a non-zero module over $R_P$ cannot exceed
  $\depth(R_P)$.
\end{proof}

\begin{Definition} \label{def42}
  A ring $R$ is $q$-Gorenstein ($q \geq 0$ an integer) if it satisfies
  Serre's condition $S_q$ and is Gorenstein in codimension $q-1$;
  equivalently, $R$ is $q$-Gorenstein if $R_P$ is Gorenstein for every prime
  ideal $P$ of $R$ with $\depth(R_P) \leq q-1$.
\end{Definition}

In particular, all rings are $0$-Gorenstein, and a ring is $1$-Gorenstein
if and only if it satisfies $S_1$ and is Gorenstein in codimension 0 (this
includes, in particular, all reduced rings).

Notice the similarity with the definition of $q$-regular rings: a ring is
$q$-regular if $R_P$ is regular whenever $\depth(R_P) \leq q-1$, or
equivalently, if $R$ satisfies $S_q$ and $R_{q-1}$ (thus 1-regular is the same
as reduced, and 2-regular is the same as normal). Of course, $q$-regular
implies $q$-Gorenstein.

\begin{Corollary} \label{cor43}
  If $R$ is $(k-1)$-Gorenstein and $M$ is an $R$-module, then $M$ is
  $k$-torsionless if and only if it is a $k^{\text{th}}$ syzygy.
\end{Corollary}

For example, over a reduced ring $R$, every second syzygy (and in particular
every dual) is reflexive.




\end{document}